# The Probabilities of Large Deviations Associated with Multinomial Distributions

Sherzod M. Mirakhmedov

V.I.Romanovskiy Institute of Mathematics. Academy of Sciences of Uzbekistan
University str.,46. Tashkent-l00174
e-mail: shmirakhmedov@yahoo.com

**Abstract**. We consider a multinomial distribution, where the number of cells increases and the cell-probabilities decreases as the number of observations grows. The probabilities of large deviations of statistics, which has form of a sum of Borel functions of cell-frequencies, are studied. Large deviation results for the power-divergence statistics and its most popular special variants, as well as for some count statistics are derived as consequences of general theorems.

**Key words and phrases**. Chi-square statistic, Count-statistics, Log-likelihood ration statistic, large deviations, multinomial distribution, Poisson distribution, power divergence statistics.

**MSC (2000):** 60F10, 62E20, 62G20

## 1. Introduction

Let $\eta = (\eta_1,...,\eta_N)$ be a random vector of frequencies of a multinomial model $M(n,N,\mathrm{P})$ on $n \geq 1$ observations classified into $N > 1$ cells with the cell probabilities $\mathrm{P} = (p_1,...,p_N)$, $p_1 + ... + p_N = 1$, all $p_j > 0$. Hence $n = \eta_1 + ... + \eta_N$ and

$$P\{\eta_1 = m_1,...,\eta_N = m_N\} = \frac{n!}{m_1!\cdot...\cdot m_N!} p_1^{m_1} \cdot...\cdot p_N^{m_N},$$

where arbitrary non-negative integer $m_j$s are such that $m_1 + ... + m_N = n$. It is convenient to describe this model in terms of a random allocation of particles into cells: $n$ particles are allocated into $N$ cells indexed 1 through $N$ at random, successively and independently of each other, the probability of a particle falls into $l$-th cell is $p_l > 0$, $l = 1,...,N$, $p_1 + ... + p_N = 1$, then $\eta_l$ is the number of particles in the $l$–th cell after allocation of all $n$ particles.

We assume that $N = N(n) \to \infty$ and all $p_m = p_m(n) \to 0$, as $n \to \infty$, and consider a general class of statistics of the form

$$R_N(\eta) = \sum_{l=1}^{N} h_l(\eta_l), \qquad (1.1)$$

where $h_1(x),...,h_N(x)$ (may depend on $N$ and $n$) are real-valued functions defined on the non-negative axis. These functions can be random, then it is assumed that for any $x_1,...,x_N$, random variables (r.v.) $h_1(x_1),...,h_N(x_N)$ are mutually independent and independent of $(\eta_1,...,\eta_N)$. We will pay special attention to the family of power-divergence statistics (PDS) of Cressie and Read (1984), which is a sub-class of (1.1) with $h_l(x) = 2x\left((x/np_l)^d - 1\right)/d(d+1)$. The most notably special versions of PDS are

$$\chi_N^2 = \sum_{m=1}^{N} \frac{(\eta_m - np_m)^2}{np_m}, \quad \Lambda_N = 2\sum_{m=1}^{N} \eta_m \ln\frac{\eta_m}{np_m}, \quad T_N^2 = 4\sum_{m=1}^{N} (\sqrt{\eta_m} - \sqrt{np_m})^2. \quad (1.2)$$

These statistics are known as Pearson's chi-squared, the log-likelihood ratio and the Freeman-Tukey statistics, respectively. The class of count statistics (CS) forms another important variant of (1.1). In particularly the following CS are of interest:

$$\mu_r = \sum_{m=1}^{N} I\{\eta_m = r\}, \quad w_r = \sum_{m=1}^{N} I\{\eta_m \geq r\}, \quad C_n = \sum_{m=1}^{N} (\eta_m - 1) I\{\eta_m > 1\}, \quad (1.3)$$

where $r = 0, 1, ..., n$ and $I\{A\}$ is the indicator function of $A$. In terms of the aforementioned random allocation scheme the statistics $\mu_r$, $w_r$, and $C_n$ counts, respectively: the number of cells that contain exactly $r \geq 0$ particles, the number of cells that contain at least $r \geq 1$ particles, and the number of collisions (i.e. the number of times a particle falls in a cell that already has a particle in it). Consider now a random allocation with random level. Let the cell with index $m$ be assigned a non-negative, integer valued r.v. $v_m$, $m = 1, ..., N$. The cell with index $m$ is said to be filled up if $\eta_m \geq v_m$ after allocation all the particles. A special version of (1.1) that of interest is the number of unfilled cells $\Phi_N$, say, viz.

$$\Phi_N = \sum_{m=1}^{N} B_m(\eta_m), \quad (1.4)$$

where $B_m(x)$ is a Bernoulli r.v.: $P\{B_m(x) = 1\} = P\{v_m > x\}$, $P\{B_m(x) = 0\} = P\{v_m \leq x\}$, $m = 1, ..., N$.

If $h_m(x) = h(x)$ for all $m$ then (1.1) is said to be symmetric statistics. Obviously $\chi_N^2$, $\Lambda_N$ and $T_N^2$ are symmetric if $p_1 = ... = p_N = N^{-1}$. All statistics (1.3) are symmetric, whereas statistic $\Phi_N$ is symmetric if $v_1, ..., v_N$ are i.i.d. r.v.s.. PDS, specifically statistics (1.2), play a leading role in goodness-of-fit tests on grouped data, whereas CS (1.3), (1.4) are used in problems associated with occupancy problems, see Kolchin et al (1976), Ivanov V.A. (1985), Khmaladze (2011). We also draw the attention of readers to the article by L'ecuyer et al (2002), where the statistics (1.1) and its special versions, including (1.2) and (1.3) are used in construction the serial tests for testing of uniformity and independence of the output sequence of general-purpose uniform random number generators.

In the present paper we are interested in $P\{\tilde{R}_N(\eta) \geq x_N\}$, where $\tilde{R}_N(\eta)$ is a standardized version of $R_N(\eta)$ and $x_N \to \infty$, as $n \to \infty$. This probability of large deviation has been studied earlier by Sirajdinov et al (1989) and Ivchenko and Mirakhmedov (1995) and Mirakhmedov (2020). The last article contains a few typos, which are corrected in Mirakhmedov's (2021). In course of revision of that paper the author observed that several refinements of the results on applications of Theorem 2.2 could be made. The purpose of this paper is to present refined results, therefore, it can be considered as updated variant of article Mirakhmedov (2020). As such, we send readers to Mirakhmedov (2020) for motivation and review of the literature on the topic.

The rest of the paper is organized as follows. The general results are presented in Section 2; these theorems are applied for the class of PDS in Section 3, and in Section 4 for the statistics (1.2), (1.3) and (1.4); the proofs are presented in Section 5. The auxiliary Assertions are collected in Appendix.

Throughout the paper we use $c_i$, $C_i$ to denote absolute constants whose value may differ at each occurrence, $\varsigma \sim F$ stands for "r.v. $\varsigma$ has the distribution $F$", $Poi(\lambda)$ - the Poisson distribution with parameter $\lambda > 0$ and $\Phi(u)$ - the standard normal distribution function. Everywhere in the next sections $\xi, \xi_1, ..., \xi_N$ are independent r.v.s, where $\xi_m \sim Poi(np_m)$; $\xi \sim Poi(\lambda_n)$, $\lambda_n = n/N$, $\pi_l(\lambda) = \lambda^l e^{-\lambda}/l!$, $p_{\max} = \max_{1 \leq k \leq N} p_k$, $p_{\min} = \min_{1 \leq k \leq N} p_k$. All asymptotic statements are considered as $n \to \infty$.

## 2. General Results. Set

$$g_m(\xi_m) = h_m(\xi_m) - Eh_m(\xi_m) - \tau_n(\xi_m - np_m),$$

$$A_N = \sum_{m=1}^{N} Eh_m(\xi_m), \quad \tau_n = n^{-1}\sum_{m=1}^{N} \text{cov}\left(h_m(\xi_m), \xi_m\right),$$

$$\tilde{\sigma}_N^2 = \sum_{m=1}^{N} Var h_m(\xi_m), \quad \sigma_N^2 = \sum_{m=1}^{N} Var g_m(\xi_m) = \tilde{\sigma}_N^2 - n\tau_n^2. \tag{2.1}$$

Note that under very general set-up $A_N$ and $\sigma_N^2$ becomes the asymptotical value of $ER_N(\eta)$ and $VarR_N(\eta)$, respectively. In particularly, under the conditions of below Theorems 2.1 and 2.2

$$ER_N(\eta) = A_N + o(N), \quad VarR_N(\eta) = \sigma_N^2(1+o(1)). \tag{2.2}$$

The following theorem was proved by Ivchenko and Mirakhmedov (1995), see also Siragdinov et al (1989).

**Theorem 2.1.** Assume that $n \to \infty$, $N \to \infty$ such that

$$n/N \to \lambda \in (0, \infty), \tag{2.3}$$

$$Np_{\max} \leq c_1, \tag{2.4}$$

$$\liminf N^{-1}\sigma_N^2 > 0, \tag{2.5}$$

$$\max_{1 \leq m \leq N} Ee^{H|h_m(\xi_m)|} \leq c_2, H > 0, c_2 > 0. \tag{2.6}$$

Then uniformly in $x_n \geq 0$, $x_n = o(\sqrt{N})$

$$P\{R_N > x_n\sigma_N + A_N\} = (1-\Phi(x_n))\exp\{M_n(x_n)\}\left(1+o\left(\frac{1+x_n}{\sqrt{N}}\right)\right), \quad (2.7)$$

$$P\{R_N < -x_n\sigma_N + A_N\} = \Phi(-x_n)\exp\{M_n(-x_n)\}\left(1+o\left(\frac{1+x_n}{\sqrt{N}}\right)\right), \quad (2.8)$$

where $M_n(u) = u^3(\mu_{0n} + \mu_{1n}u + ...)$ is a series where $|\mu_{jn}| \leq \mu_j N^{-(j+1)/2}$ for sufficiently large $N$, $\mu_j$ do not depend on $N$. In particularly, due to Sirajdinov et al (1989)

$$\mu_{0n} = \frac{\beta_{3N}}{6\sigma_N^3}, \quad \mu_{1n} = \frac{\beta_{4N}}{24\sigma_N^4} - \frac{\beta_{3N}^2}{8\sigma_N^6} + \frac{1}{n\sigma_N^4}\left(\sum_{m=1}^{N} Eg_m^2(\xi_m)(\xi_m - np_m)\right)^2 - \frac{1}{8\sigma_N^4}\sum_{m=1}^{N}(Eg_m^2(\xi_m))^2.$$

where $\beta_{l,N} = Eg_1^l(\xi_1) + ... + Eg_N^l(\xi_N)$.

Set $K_n(a,b) = \left(n^{1-b}p_{\max}^{-b}\right)^{1/(1+\bar{a})}$, where integers $a \geq 0$, $b \geq 0$ and $\bar{a} = \max(1,a)$.

**Theorem 2.2.** Let $N \to \infty$, $p_{\max} \to 0$ as $n \to \infty$, and the functions $h_m(\cdot)$ be non-negative. If

(i) For each integer $s \in [3, k_n]$ there exist non-negative $a_1$, $a_2$ and $b_1$, $b_2$ such that

$$E\left((\xi_m - np_m)^2 h_m^s(\xi_m)\right) = O\left(s^a (np_m)^b Eh_m^s(\xi_m)\right), \quad (2.9)$$

where $a = a_1$, $b = b_1$ for all $m \in \mathcal{N} \subseteq (1,2,...,N)$, and $a = a_2$, $b = b_2$ for all $m \in (1,2,...,N)\setminus\mathcal{N}$, and

$$k_n = o\left(\min(p_{\max}^{-1}, K_n(a_1,b_1), K_n(a_2,b_2))\right), \quad (2.10)$$

(ii) For each integer $s \in [3, k_n]$ there exists a sequence of positive numbers $V_n$ such that for some $\nu \geq 0$

$$|\mathcal{C}_s(h_m(\xi_m))| \leq (s!)^{1+\nu} V_n^{s-2} Var h_m(\xi_m), \quad m = 1,...,N. \quad (2.11)$$

Then for all $x_n$ such that

$$0 \leq x_n = o(\Upsilon_n), \quad (2.12)$$

where $\Upsilon_n = \min\left(W_n^{1/(1+2\nu)}, k_n^{1/2}\right)$, $W_n = \sigma_N^3 / V_n \tilde{\sigma}_N^2$, it holds

$$P\{R_N(\eta) > x_n\sigma_N + A_N\} = (1-\Phi(x_n))\exp\{M_{\nu,n}(x_n)\}\left(1+O\left(\frac{x_n+1}{\Upsilon_n}\right)\right), \quad (2.13)$$

and

$$P\{R_N(\eta) < -x_n\sigma_N + A_N\} = \Phi(-x_n)\exp\{M_{\nu,n}(-x_n)\}\left(1+O\left(\frac{x_n+1}{\Upsilon_n}\right)\right), \quad (2.14)$$

where $M_{\nu,n}(x) = 0$ if $\nu \geq 1$, $M_{\nu,n}(x) = x^3(\mu_{0n} + ... + \mu_{ln}x^l)$, where $0 \leq l < \kappa(\nu)$, $\kappa(0) = \infty$, $\kappa(\nu) = \nu^{-1} - 1$ if $\nu \in (0,1)$. Yet $M_{0,n}(x) = M_n(x)$, and it holds

$$|M_{v,n}(x)| \leq c(v)|x|^3 / W_n^{1/(1+2v)}. \tag{2.15}$$

To compare these two theorems we make the following comments. The deviation zone $x_N = o(\sqrt{N})$ of Theorem 2.1 is maximally possible for this Cramèr type result; application, however, of Theorem 2.1 is limited by sparse multinomial models only, since conditions (2.3), (2.4), and by class of statistics satisfying Cramèr condition (2.6). For instance, class of PDS with parameter $d > 0$ does not satisfy condition (2.6). In contrast, Theorem 2.2 is applicable for the arbitrary multinomial models and for statistics that do not necessarily satisfy Cramér condition. However, condition (i) turns out to be restrictive, and may narrow the deviation zone so that the exponential factor does not effect, see for instance, the applications below.

In what follows, in applications of Theorems 2.1 and 2.2, we will present the statements, which follows from (2.7) and (2.13) only, while the results similar to (2.14) hold as well.

### 3. Application to the Power Divergence Statistics (PDS).

The PDS of Cressie and Read (1984) is defined as

$$CR_N(d) = \frac{2}{d(d+1)} \sum_{l=1}^{N} \eta_l \left[ (\eta_l / np_l)^d - 1 \right] = \frac{2}{d(d+1)} \sum_{l=1}^{N} np_l (\eta_l / np_l)^{d+1} - \frac{2n}{d(d+1)}.$$

Note that for $d = 0$ the PDS is defined by continuity: $CR_N(0) = \lim_{d \to 0} CR_N(d) = \Lambda_N$, also $CR_N(1) = \chi_N^2$ and $CR_N(-1/2) = T_N^2$. We assume that $d > -1$ in order Theorems 2.1 and 2.2 to be applicable. Further, since we deal with standardized variant of PDS, we can, whenever it is necessary, consider PDS as a special version of statistics $R_N(\eta)$ with kernel functions $h_l(x) = h_{d,l}(x)$, where

$$h_{d,l}(x) = np_l (x / np_l)^{d+1}, d > -1, d \neq 0, \quad \text{else } h_l(x) = h_{0,l}(x) = 2x \log(x / np_l). \tag{3.1}$$

The notation (2.1) for this case has the following form: for $d > -1$, $d \neq 0$:

$$A_N(d) = \sum_{m=1}^{N} (np_m)^{-d} E\xi_m^{d+1}, \tag{3.2}$$

$$\tau_n(d) = \frac{1}{n} \sum_{m=1}^{N} (np_m)^{-d} E\xi_m^{d+1}(\xi_m - np_m), \tag{3.3}$$

$$\tilde{\sigma}_N^2(d) = \sum_{m=1}^{N} (np_m)^{-2d} Var \xi_m^{d+1}. \tag{3.4}$$

and

$$A_N(0) = 2n \sum_{m=1}^{N} p_m E\left[ (\xi_m / np_m) \ln(\xi_m / np_m) \right], \tag{3.5}$$

$$\tau_n(0) = 2n^{-1} \sum_{m=1}^{N} E\left( (\xi_m - np_m) \xi_m \ln(\xi_m / np_m) \right), \tag{3.6}$$

$$\tilde{\sigma}_N^2(0) = 4n^2 \sum_{m=1}^{N} p_m^2 Var\left[(\xi_m / np_m)\ln(\xi_m / np_m)\right]. \tag{3.7}$$

Again we remind that $A_N(d)$ and

$$\sigma_N^2(d) = \tilde{\sigma}_N^2(d) - n\tau_n^2(d) \tag{3.8}$$

are asymptotic expectation and asymptotic variance of $CR_N(d)$, respectively.

**3.1. The sparse multinomial model**: $n/N \to \lambda \in (0, \infty)$. Theorems 2.1 and 2.2 imply

**Theorem 3.1**. Suppose $n \to \infty$, $N \to \infty$ such that

$$n/N \to \lambda \in (0, \infty) \text{ and } c_3 \leq Np_m \leq c_4, \text{ some } c_3 > 0, c_4 > 0. \tag{3.9}$$

(i) Let $-1 < d \leq 0$. Then uniformly in $x_n \geq 0$, $x_n = o(\sqrt{N})$ for the $P\{CR_N(d) > x_n \sigma_N(d) + A_N(d)\}$ and $P\{CR_N(d) < -x_n \sigma_N(d) + A_N(d)\}$ the relations (2.7) and (2.8), respectively, hold.

(ii) Let $d > 0$. Then uniformly in $x_n \geq 0$, $x_n = o(\min(N^{1/6}, N^{1/2(1+2d)}))$ one has

$$P\{CR_N(d) > x_n \sigma_N(d) + A_N(d)\} = (1 - \Phi(x_n))(1 + o(1)), \tag{3.10}$$

**Corollary 3.1**. Assume (3.9) is fulfilled. If $-1 < d \leq 0$ then for $x_n \to \infty$, $x_n = o(\sqrt{N})$, and if $d > 0$ then for $x_n \to \infty$, $x_n = o(\min(N^{1/6}, N^{1/2(1+2d)}))$ one has

$$\log P\{CR_N(d) > x_n \sigma_N(d) + A_N(d)\} = -\frac{x_n^2}{2}(1 + o(1)). \tag{3.11}$$

**3.2. The dense multinomial model**: $\lambda_n \to \infty$. In this case, each PDS asymptotically in distribution coincides with the chi-square statistics.

**Theorem 3.2**. Let $p_{\max} = o(1)$ and $np_{\min} \to \infty$. Then for each $d > -1$ uniformly in $x_n \geq 0$, $x_n = o\left(\min(N^{1/6}, p_{\max}^{-1/4})\right)$ one has

$$P\{CR_N(d) > x_n\sqrt{2N} + N\} = (1 - \Phi(x_n))(1 + o(1)), \tag{3.12}$$

**Corollary 3.2**. Let $p_{\max} = o(1)$ and $np_{\min} \to \infty$. Then for each $d > -1$ and $x_n \to \infty$, $x_n = o(\min(N^{1/6}, p_{\max}^{-1/4}))$ one has (3.12).

**Remark 3.1**. From $np_{\min} \to \infty$ it follows $\lambda_n \to \infty$. The condition $np_{\min} \to \infty$ can be replaced by $\lambda_n \to \infty$ and $Np_{\min} \geq c > 0$, which imply in turn $np_{\min} \to \infty$. It is obvious that if $p_{\max} = O(N^{-2/3})$, then relations (3.12) is hold uniformly in non-negative $x_N = o(N^{1/6})$.

**3.3. The very sparse multinomial model**: $\lambda_n \to 0$. We assume that $Np_{\max} \leq C$, which in our case, yields $np_{\max} \to 0$. The $CR_N(d)$, being divergence measure between proposed and observed models, depends on the intended model. As it turns out, this fact plays a role for very sparse models; for example, the asymptotic value of the expectation and the variance of PDS are different for non-uniform and uniform multinomial models. The standardized version of $CR_N(d)$ coincides with standardized version of $R_N^d = h_{d,1}(\eta_1) + ... + h_{d,N}(\eta_N)$, where $h_{d,m}(x)$ is defined in (3.1). Set $P_{jN}(a) = p_1^{j-a} + ... + p_N^{j-a}$. The asymptotic value of the expectation and variance of $R_N^d$, when $d \neq 0$, are equal, respectively to

$$A_N(d) = n^{1-d} P_{1N}(d), \tag{3.13}$$

$$\tilde{\sigma}_N^2(d) = n^{1-2d} P_{1N}(2d), \tag{3.14}$$

$$\sigma_N^2(d) = n^{1-2d}\left[P_{1N}(2d) - P_{1N}^2(d)\right]$$
$$+ 2n^{2(1-d)}\left[(2^{2d} - 1)P_{2N}(2d) - 2(2^d - 1)P_{1N}(d)P_{2N}(d)\right], \ d > -1, d \neq 0. \tag{3.15}$$

**Remark 3.2.** In the formula for $\sigma_N^2(d)$ given in Mirakhmedov (2020, Eq. (2.6)), the first term of (3.14) is missed. Let $P\{Y = p_m^{-d}\} = p_m, m = 1,...,N$. Then $P_{1N}(2d) - P_{1N}^2(d) = VarY$, and hence $P_{1N}(2d) - P_{1N}^2(d) = 0$ if and only if $(p_1,..., p_N) = (N^{-1},..., N^{-1})$. Next, $nP_{2N}(2d) \leq np_{\max} P_{1N}(2d)$ and $nP_{1N}(d)P_{2N}(d) \leq np_{\max} P_{1N}^2(d)$.

Consider now $R_N^0 = \Lambda_N$. Let $(p_1,..., p_N) \neq (N^{-1},..., N^{-1})$, a r.v. $Z_N$ be such that $P\{Z_N = -\ln np_m\} = p_m, m = 1,...,N$. For the asymptotic value of the expectation and variance of $R_N^0$ we have

$$A_N(0) = 2nEZ_N, \ \tilde{\sigma}_N^2(0) = 4nEZ^2, \ \sigma_N^2(0) = 4n\,\text{var}\,Z_N. \tag{3.16}$$

Let now, $(p_1,..., p_N) = (N^{-1},..., N^{-1})$. Then

$$A_N(d) = n\lambda_n^{-d}, \ \sigma_N^2(d) = 2(2^d - 1)^2 n\lambda_n^{1-2d}, \ d > -1, \ d \neq 0, \tag{3.17}$$

$$A_N(0) = 2\ln 2\, n\lambda_n, \ \tilde{\sigma}_N^2(0) = \sigma_N^2(0) = 8\ln^2 2\, n\lambda_n. \tag{3.18}$$

Set $\bar{R}_N^d = \eta_1^{1+d} + ... + \eta_N^{1+d}$, $d > -1$, $d \neq 0$, else $\bar{R}_N^0 = 2\eta_1 \ln \eta_1 + ... + 2\eta_N \ln \eta_N$. For the uniform model standardized version of $R_N^d$ and $\bar{R}_N^d$ coincides. For $\bar{R}_N^d$ the quantities (3.2) - (3.7) are:

$$A_N(d) = n, \ \tilde{\sigma}_N^2(d) = n, \ \sigma_N^2(d) = 2(2^d - 1)^2 n\lambda_n, \ d \neq 0, \tag{3.19}$$

else $A_N(0)$ and $\sigma_N^2(0)$ are as in (3.18). Put $\Upsilon_n(d) = \min(n^{1/4}, (W_N(d))^{1/(1+2d^*)})$, where $d^* = \max(0, d)$, $W_N(d) = \sigma_N^3(d)(np_{\min})^d (\tilde{\sigma}_N^2(d))^{-1}$ if $d \neq 0$, else $W_N(0) = \sigma_N^3(0)(\tilde{\sigma}_N^2(0)|\ln np_{\min}|)^{-1}$.

**Theorem 3.3.** Let $(p_1,..., p_N) \neq (N^{-1},..., N^{-1})$, $\lambda_n \to 0$, $n\lambda_n \to \infty$ and $np_{\max} \to 0$ (or $Np_{\max} \leq C$). Then uniformly in $x_n \geq 0$, $x_n = o(\Upsilon_n(d))$ it holds

$$P\{R_N^d > x_n \sigma_N(d) + A_N(d)\} = (1-\Phi(x_n))\exp\{M_{d^*,n}(x_n)\}\left(1+O\left(\frac{x_n+1}{\Upsilon_n(d)}\right)\right).$$

Set $\bar{\Upsilon}_n(d) = \min(n^{1/4}, (n\lambda_n^3)^{1/2(1+2d^*)})$.

**Theorem 3.4.** Let $p_1 = ... = p_N = N^{-1}$, $\lambda_n \to 0$ and $n\lambda_n^3 \to \infty$. Then uniformly in $x_n \geq 0$, $x_n = o(\bar{\Upsilon}_n(d))$ it holds

(i) $P\{\bar{R}_N^d > x_n |2^d - 1|\sqrt{2n\lambda_n} + n\} = (1-\Phi(x_n))\exp\{M_{d^*,n}(x_n)\}\left(1+O\left(\frac{x_n+1}{\bar{\Upsilon}_n(d))}\right)\right), d \neq 0$ ,

(ii) $P\{\bar{R}_N^0 > x_n 2\sqrt{2\ln 2n\lambda_n} + 2\ln 2\, n\lambda_n\} = (1-\Phi(x_n))\exp\{M_{0,n}(x_n)\}\left(1+O\left(\frac{x_n+1}{\bar{\Upsilon}_n(0)}\right)\right).$

**Remark 3.3.** Let $p_m = (1-\alpha)/N^{1-\alpha}m^\alpha$, $m=1,...,N$ and $\alpha \in [0,1)$. This example of the multinomial model has been considered by Rempala and Wesolowski (2016), as it is of interest when testing for signal-noise threshold in data with large number of support points (Pietrzak et al., 2016). We have $np_{\min} = (1-\alpha)\lambda_n$, $p_{\max} = (1-\alpha)/N^{1-\alpha}$, $Np_m = (1-\alpha)N^\alpha m^{-\alpha}$, $1-\alpha \leq Np_m \leq (1-\alpha)N^\alpha$. Theorem 3.2 for this model states that if $\lambda_n \to \infty$ then the relation (3.12) hold uniformly in $x_n \geq 0$ such that: $x_n = o(N^{1/6})$ if $\alpha \in [0,1/3]$, and $x_n = o(N^{(1-\alpha)/4})$ if $\alpha \in (1/3,1)$. Theorem 3.3 can be applied if $\alpha \in [0,1/2)$. However, Theorem 3.1 possible to apply for $\alpha = 0$, i.e. for uniform model only.

**Corollary 3.3.** Let $p_1 = ... = p_N = N^{-1}$, $\lambda_n \to 0$ and $n\lambda_n^3 \to \infty$. Then uniformly in $x_n \geq 0$, $x_n = o(\tilde{W}_N(d))$ it holds

$$\log P\{\bar{R}_N^d > x_n |2^d - 1|\sqrt{2n\lambda_n} + n\} = -\frac{x_n^2}{2}(1+o(1)), d \neq 0,$$

$$\log P\{\bar{R}_N^0 > x_n 2\sqrt{2\ln 2n\lambda_n} + 2\ln 2\, n\lambda_n\} = -\frac{x_n^2}{2}(1+o(1)).$$

We end Section 3 by noting that the theorems of this section should be regarded as results under the null hypothesis that the cell probabilities are $p_1,...,p_N$.

**4. Application to the special variants of the statistic $R_N(\eta)$.**

Consider application of results of Sections 2 and 3 to the statistics (1.2), (1.3) and (1.4). Remind that (1.2) statistics are PDS: $\chi_N^2 = CR_N(1)$, $\Lambda_N = CR_N(0)$ and $T_N^2 = CD_N(-1/2)$, whereas the (1.3) and (1.4) statistics are out of the class of PDS.

**4.1. Chi-square statistic.** For the $\chi_N^2$ statistic $h_m(u) = (u - np_m)^2 / np_m$,

$$\tilde{\sigma}_N^2 = \sum_{m=1}^{N} \frac{1}{np_m} + 2N, \quad \sigma_N^2 = 2N + \sum_{m=1}^{N}\left(\frac{1}{np_m} - \frac{1}{\lambda_n}\right) = \tilde{\sigma}_N^2 - N\lambda_n^{-1}. \tag{4.1}$$

As consequence of Theorem 3.1 (ii) and Theorem 3.2 we obtain

**Theorem 4.1.** Let $n \to \infty$, $N \to \infty$ such that the condition (3.9) or

$$\lambda_n \to \infty \text{ and } c_1 \le Np_{\min} \le Np_{\max} \le c_2 N^{1/3}, \tag{4.2}$$

is fulfilled. Then for all $x_n \ge 0 : x_n = o(N^{1/6})$ one has

$$P\{\chi_N^2 > x_n \sigma_N + N\} = (1 - \Phi(x_n))(1 + o(1)). \tag{4.3}$$

**Remark 4.1.** Let $Y$ be a r.v. such that $P\{Y = p_m^{-1}\} = p_m$, $m = 1, ..., N$. Then $EY = N$ and

$VarY = \sum_{m=1}^{N} p_m^{-1} - N^2 = n\sum_{m=1}^{N}\left((np_m)^{-1} - \lambda_n^{-1}\right) \ge 0$. Hence $\sigma_N^2 \ge 2N$. Under the conditions (4.2)

$np_{\min} \to \infty$, and hence $\sigma_N^2 = 2N(1 + o(1))$.

**Theorem 4.2.** Let $(p_1, ..., p_N) \ne (N^{-1}, ..., N^{-1})$ and $\lambda_n \to 0$, $n\lambda_n \to \infty$ and $Np_{\max} \le C$. Then for all $x_n \ge 0$, $x_n = o(\lambda_n \sigma_N^3 / \tilde{\sigma}_N^2)$ the Eq. (4.3) hold.

**Theorem 4.3.** Let $p_1 = ... = p_N = N^{-1}$, $\lambda_n \to 0$ and $n\lambda_n^3 \to \infty$. Then uniformly in $x_n \ge 0$, $x_n = o\left((n\lambda_n^3)^{1/6}\right)$ it holds

$$P\{\chi_N^2 > x_n \sqrt{2N} + N\} = (1 - \Phi(x_n))(1 + o(1)).$$

Theorem 2.2 allows deriving relations (4.3) under weaker conditions than Theorem 4.1. Set $\nabla_n = \max\left(1, (np_{\min})^{-1}\right)$.

**Theorem 4.4.** Let $p_{\max} = o(1)$. Then it holds Eq. (4.3) for all $x_n \ge 0$,

$$x_n = o\left(\min\left(\left(\sigma_N^3 / \tilde{\sigma}_N^2 \nabla_n\right)^{1/3}, n^{1/6}, p_{\max}^{-1/4}\right)\right).$$

**Corollary 4.1.** Let $np_{\min} \ge c_0 > 0$. Then Eq. (4.3) hold uniformly in $x_n \ge 0$, $x_n = o(\min(N^{1/6}, p_{\max}^{-1/4}))$, and if additionally $p_{\max} = O(N^{-2/3})$ then for all $x_n \ge 0$, $x_N = o(N^{1/6})$.

Let's consider the multinomial model considered in Remark 3.2. Corollary 4.1 says that if $\lambda_n \ge c_0$, some $c_0 > 0$ (the condition weaker than $\lambda_n \to \infty$ as in Remark 3.2) then Eq. (4.3) is valid for all $x_n \ge 0$, $x_n = o(\min(N^{1/6}, N^{(1-\alpha)/4}))$.

The following consequence of Theorem 4.4 is useful in studying the intermediate asymptotic efficiency of the chi-squared test in verifying uniformity of a discrete distribution versus a family of sequences of alternatives approaching the hypothesis, see Appendix A2.

**Corollary 4.2.** Let $p_m = N^{-1}\left(1 + \delta(n)\ell_{m,n}\right)$, $m = 1, ..., N$, where $\delta(n) \to 0$ and

$$\sum_{m=1}^{N} \ell_{m,n} = 0, \quad \frac{1}{N}\sum_{m=1}^{N} \ell_{m,n}^2 = \ell^2 < \infty.$$

Then for arbitrary $\lambda_n$ and $x_n \to \infty$, $x_n = o\left(N^{1/6}\min(1, \lambda_n^{2/3})\right)$ it holds

$$\log P\{\chi_N^2 > x_n\sqrt{2N} + N\} = -\frac{x_n^2}{2}(1+o(1)).$$

**Remark 4.2**. From (4.3) when $x_n \to \infty$ we obtain $\log P\{\chi_N^2 > x_n\sigma_N + N\} = -x_n^2/2(1+o(1))$. From Kallenberg (1985, Eq. (2.17)) this relation follows if $\log N = o(x_n^2)$, $x_n = o(N^{1/6})$ and $N = o(n^{3/8})$ (i.e. the multinomial model is so dense that $\lambda_n \gg n^{5/8}$). Last condition excludes, for instance, the case recommended by Mann and Wald (1942), who obtained the relation $N = cn^{2/5}$ concerning the optimal choice of the number of groups in chi-square goodness of fit test.

**4.2. Likelihood ratio statistic $\Lambda_N$.** Remind that $\Lambda_N = R_N^0$, and hence Theorems 3.3 and 3.4 consist results for statistic $\Lambda_N$ in very sparse models. Next, Theorem 3.1 and 3.2 imply

**Theorem 4.5.** (i) Let condition (3.9) is fulfilled. Then for all $x_n \geq 0$, $x_n = o(\sqrt{N})$ it holds

$$P\{\Lambda_N > x_n\sigma_N(0) + A_N(0)\} = (1 - \Phi(x_n))\exp\{M_{0,n}(x_n)\}\left(1 + O\left(\frac{1+x_n}{\sqrt{N}}\right)\right).$$

(ii) Let $p_{\max} = O(N^{-2/3})$ and $np_{\min} \to \infty$. Then for all $x_n \geq 0$, $x_n = o(N^{1/6})$ one has

$$P\{\Lambda_N > x_N\sqrt{2N} + N\} = (1 - \Phi(x_N))(1+o(1)). \tag{4.4}$$

**Corollary 4.3**. If condition (3.9) is fulfilled then for $x_n \to \infty$, $x_n = o(\sqrt{N})$, and if $np_{\min} \to \infty$ and $p_{\max} = O(N^{-2/3})$ then for $x_n \to \infty$, $x_N = o(N^{1/6})$ it holds

$$\log P\{\Lambda_N > x_n\sqrt{2N} + N\} = -\frac{x_n^2}{2}(1+o(1)).$$

**Remark 4.3.** Theorem 4.5 (ii) imply $\log P\{\Lambda_N > x_n\sqrt{2N} + N\} == -x_n^2/2(1+o(1))$ if $x_N \to \infty$ and $x_n = o(N^{1/6})$. From Kallenberg (1985, Eq. (2.13)) such relation follows under additionally conditions $N = o(n^{3/7})$ and $\log N = o(x_n^2)$.

**4.3. The Freeman-Tukey statistic $T_N^2$** is $CR_N(-1/2)$ statistic, therefore directly from Theorems 3.1 and 3.2 it follows

**Theorem 4.6**. (i) Let condition (3.9) is fulfilled. Then for the $P\{T_N^2 > x_n\sigma_N(-1/2) + A_N(-1/2)\}$ it holds Eq. (2.7) uniformly in $x_n \geq 0$, $x_n = o(\sqrt{N})$.

(ii) Let $p_{\max} = O(N^{-2/3})$ and $np_{\min} \to \infty$. Then uniformly in $x_n \geq 0$, $x_n = o(N^{1/6})$ it holds

$$P\{T_N^2 > x_n\sqrt{2N} + N\} = (1-\Phi(x_n))(1+o(1)).$$

(iii) Let $(p_1,...,p_N) \neq (N^{-1},...,N^{-1})$, $\lambda_n \to 0$, $n\lambda_n \to \infty$ and $Np_{\max} \leq C$. Then

$$P\left\{\frac{T_N^2 - 8(n - A_N(-1/2))}{8\sigma_N(-1/2)} > x_n\right\} = (1-\Phi(x_n))\exp\{M_{0,n}(x_n)\}\left(1+O\left(\frac{x_n+1}{\Upsilon_n(-1/2)}\right)\right),$$

uniformly in $x_n \geq 0$, $x_n = o(\Upsilon_n(-1/2))$, where $\Upsilon_n(-1/2) = \min(n^{1/4}, W_n(-1/2))$, $A_N(-1/2)$, $\sigma_N^2(-1/2)$ and $\tilde{\sigma}_N^2(-1/2)$ are defined in (3.15), (3.16) and (3.17), respectively, and $W_n(-1/2) = \sigma_N^3(-1/2)(\sqrt{np_{\min}}\tilde{\sigma}_N^2(-1/2))^{-1}$.

(iv) Let $p_1 = ... = p_N = N^{-1}$, $\lambda_n \to 0$ and $n\lambda_n^3 \to \infty$. Then uniformly in $x_n \geq 0$, $x_n = o(W_n)$ it holds

$$P\left\{\frac{T_N^2 - 8n(1-\sqrt{\lambda_n})}{8(\sqrt{2}-1)\sqrt{n}\lambda_n^{3/2}} > x_n\right\} = (1-\Phi(x_n))\exp\{M_{0,n}(x_n)\}\left(1+O\left(\frac{x_n+1}{W_n}\right)\right),$$

where $W_n = n^{1/4}$ if $N \geq n^{7/6}$ and $W_n = \sqrt{n\lambda_n^3}$ if $N \leq n^{7/6}$.

**4.5. The count statistics.**

**4.5.1.** *First we consider* CS (1.3). Remind, for the $\mu_r$ statistic $h_m(x) = I\{x = r\}$. Set $\pi_s(\lambda) = \lambda^s e^{-\lambda}/s!$. The notation (2.1) in this case has the following form:

$$A_N = A_{rN} = \sum_{m=1}^{N}\pi_r(np_m), \quad \tau_N = \tau_{rN} = n^{-1}\sum_{m=1}^{N}(r-np_m)\pi_r(np_m),$$

$$\tilde{\sigma}_N^2 = \tilde{\sigma}_{rN}^2 = \sum_{m=1}^{N}\pi_r(np_m)(1-\pi_r(np_m)). \quad \sigma_N^2 = \sigma_{rN}^2 = \tilde{\sigma}_{rN}^2 - n\gamma_{rN}^2.$$

Set $\Upsilon_{r,n} = \min(\sigma_{rN}^3/\tilde{\sigma}_{rN}^2, n^{1/4}, (np_{\max}^2)^{-1/4})$.

**Theorem 4.7.** Let fixed $r \geq 0$. Uniformly in $x_n \geq 0$, $x_n = o(\Upsilon_{r,n})$ one has

$$P\{\mu_r > x_n\sigma_{rN} + A_{rN}\} = (1-\Phi(x_n))\exp\{M_{0,n}(x)\}\left(1+o\left(\frac{x_n+1}{\Upsilon_{r,n}}\right)\right). \tag{4.5}$$

Assume that there exist $c_1 > 0$ and $c_2 > 0$ such that

$$c_1 \leq Np_{\min} \leq Np_{\max} \leq c_2. \tag{4.6}$$

Condition $\sigma_{rN} \to \infty$, that is $Var\mu_r \to \infty$, is necessary for the above results of large deviations to be valid for $\mu_r$. In turn, this condition is fulfilled for the following cases of behavior of parameter $\lambda_n$:

(i) $\lambda_n \to \lambda \in (0,\infty)$.

(ii) $\lambda_n \to 0$ and $N\lambda_n^{\tilde{r}} \to \infty$, where $\tilde{r} = \max(2,r)$.

(iii) $\lambda_n \to \infty$ but $\lambda_n \leq \delta c_2^{-1} \log N + r \log \log N$, some $\delta \in (0,1)$, $c_2$ from the condition (4.6).

For all these three cases $\sigma_{rN}^2 \sim Var\mu_r \sim E\mu_r \sim A_{rN} \to \infty$. Moreover,

in (i) case: $\tilde{\sigma}_{rN}^2 \sim \sigma_{rN}^2 \sim cN$, (4.7)

in case (ii): $\sigma_{rN}^2 \sim cN\lambda_n^{\tilde{r}}$, $\tilde{\sigma}_{rN}^2 \sim N\lambda_n^r / r!, r \geq 2$, and $\tilde{\sigma}_{rN}^2 \sim n$, $r = 0,1$, (4.8)

in case (iii) $\tilde{\sigma}_{rN}^2 \sim \sigma_{rN}^2 \sim cN^{1-\delta}$. (4.9)

These (i), (ii) and (iii) cases correspond to the central, left intermediate and right intermediate zones, respectively, defined by Kolchin et al (1976). Applying Theorem 2.1 for (i) case and Theorem 4.5 for (ii) and (iii) cases we obtain the following

**Corollary 4.3.** Let condition (4.6) is fulfilled. Then Theorem 4.7 holds, where $\Upsilon_{r,n} = \sqrt{N}$ in case (i), $\Upsilon_{r,n} = n^{1/4}$ in case (ii), $\Upsilon_{r,n} = \min(N^{(1-\delta)/2}, (N/\lambda_n)^{1/4})$ in case (iii).

**Remark 4.4.** We have $w_1 = N - \mu_0$ and $C_n = \mu_0 - (n - N)$. So $P\{C_n - EC_n > x\sqrt{VarC_n}\}$ $= P\{\mu_0 - E\mu_0 > x\sqrt{Var\mu_0}\}$ and $P\{w_1 - Ew_1 > x\sqrt{Varw_1}\} = P\{\mu_0 - E\mu_0 < -x\sqrt{Var\mu_0}\}$. Hence, assertions of Corollary 4.3 corresponding to the case $r = 0$ can be reformulated for $w_1$ and $C_n$ statistics as well.

**4.5.2.** *Statistic* (1.4) $\Phi_N$. We will consider the case when $p_1 = ... = p_N = N^{-1}$, the classical random allocation scheme, and assume that the levels $\nu_1, ..., \nu_N$ are independent copies of a non-negative integer-valued r.v. $\nu$ such that $P\{\nu = 0\} = \beta_0 < 1$. Set

$$\vartheta = \min\{l : l > 0, P\{\nu = l\}\}, \alpha(\vartheta) = P\{\nu = \vartheta\}/\vartheta!, \tau(\lambda) = \sum_{l=0}^{\infty} \pi_l(\lambda)P\{\nu > l\}.$$

Then

$$A_N = N\tau(\lambda_n), \tau_n = \tau'(\lambda_n),$$

$$\sigma_N^2 = N\left[\tau(\lambda_n)(1-\tau(\lambda_n)) - \lambda_n(\tau'(\lambda_n))^2\right], \tilde{\sigma}_N^2 = N\tau(\lambda_n)(1-\tau(\lambda_n)). \quad (4.10)$$

Here $\tau'(\lambda)$ is the derivative of $\tau(\lambda)$. It is easy to see that $\tau'(\lambda) = -\sum_{l=0}^{\infty} \pi_l(\lambda)P\{\nu = l+1\}$, and hence function $\tau(\lambda)$ is strongly decreasing if $\beta_0 < 1$. Notice that, if $\lambda_n \to 0$, then

$$\tau(\lambda_n) = 1 - \beta_0 - \alpha(\vartheta)\lambda_n^\vartheta(1 + O(\lambda_n)), \quad (4.11)$$

and hence $A_N = N(1-\beta_0)(1+O(\lambda_n^\vartheta))$, $\sigma_N^2 = N\beta_0(1-\beta_0)(1+O(\lambda_n^\vartheta))$.

In case $\lambda_n \to \infty$ the following two variants of $\tau(\lambda_n)$'s behavior are common

$$\tau(\lambda_n) = \rho\lambda_n^\chi e^{-\mu\lambda_n}\left(1 + O(\lambda_n^{-\delta})\right), \rho, \mu, \delta > 0 \text{ and } \chi \geq 0, \quad (4.12)$$

$$\tau(\lambda_n) = g\lambda_n^{-\varepsilon}\left(1+O(\lambda_n^m e^{-\lambda_n})\right), \text{ where } g, m > 0 \text{ and } \varepsilon \geq 1. \tag{4.13}$$

From Theorems 2.1 and 2.2 it follows

**Theorem 4.8.** Suppose $p_1 = ... = p_N = N^{-1}$ and the levels $v_1,...,v_N$ are i.i.d. r.v.s.. For all $x_n \geq 0$, $x_n = o(\Upsilon_n)$ it holds

$$P\{\Phi_N > x_n \sigma_N + A_N\} = (1-\Phi(x_n))\exp\{M_{0,n}(x_n)\}\left(1+O\left(\frac{x_n+1}{\Upsilon_n}\right)\right),$$

where (i) if $\lambda_n \to \lambda \in (0,\infty)$ then $\Upsilon_n = \sqrt{N}$, (ii) if $\lambda_n \to 0$ then $\Upsilon_n = n^{1/4}$, and (iii) if $\lambda_n \to \infty$ then $\Upsilon_n = \min((N/\lambda_n)^{1/4}, N^{(1-\delta)/2})$ if $\lambda_n \leq \delta\mu^{-1}\ln N + \chi \ln\ln N$ and (4.12) hold, and

$\Upsilon_n = \min((N/\lambda_n)^{1/4}, \sqrt{N\lambda_n^{-\varepsilon}})$ if (4.13) hold.

## 5. Proofs.

We still use the notation of previous sections. Additionally let $\mathcal{C}_k(\zeta)$ stands for the cumulant of $k$ th order of the r.v. $\zeta$ and

$$T_N = h_1(\xi_1) + ... + h_N(\xi_N), \quad S_N = (\xi_1 - np_1) + ... + (\xi_N - np_N);$$

also to keep notation simple we write $R_N$ instead of $R_N(\eta)$.

**Proposition 5.1.** Let the functions $h_m(\cdot)$ be non-negative and condition (i) of Theorem 2.2 be fulfilled. Then for each integer $s \in [3, k_n]$, where $k_n$ satisfies condition (2.10), it holds

$$\mathcal{C}_s(R_N) = \mathcal{C}_s(T_N)(1+o(1)). \tag{5.1}$$

**Proof.** We start from the following

**Lemma 5.1.** For any non-negative $s$ it holds

$$ER_N^s = \upsilon_n \int_{-\pi\sqrt{n}}^{\pi\sqrt{n}} ET_N^s \exp\left\{i\tau \frac{S_N}{\sqrt{n}}\right\} d\tau,$$

where

$$\upsilon_n \stackrel{def}{=} \left(2\pi\sqrt{n}P\{S_N = 0\}\right)^{-1} = \frac{n!e^n}{2\pi n^n \sqrt{n}} = \frac{1}{\sqrt{2\pi}}\left(1+o\left(\frac{1}{n}\right)\right). \tag{5.2}$$

**Proof.** It is well known that $\mathcal{L}((\eta_1,...,\eta_N)) = \mathcal{L}((\xi_1,...,\xi_N)/S_N = 0)$, where $\mathcal{L}(X)$ stands for the distribution of a random vector $X$. Hence $ER_N^k = E(T_N^k|S_N = 0)$. On the other hand $E(T_N^k e^{i\tau S_N}) = E\{e^{i\tau S_N} E(T_N^k|S_N)\}$. Now Lemma 5.1 follows by Fourier inversion. The Eq. (5.2) follows because $S_N + n \sim Poi(n)$ and Stirling's formula.

Let integer $s \in [3, k_n]$. Use Lemma 5.1 to write

$$ER_N^s = \upsilon_n \sum_{l=1}^{s} \sideset{}{'}\sum_{l,s} \sideset{}{''}\sum_{l} \int_{-\pi\sqrt{n}}^{\pi\sqrt{n}} E\left( \left(h_{j_1}(\xi_{j_1})\right)^{s_1} \cdot \ldots \cdot \left(h_{j_l}(\xi_{j_l})\right)^{s_l} \exp\left\{i\tau \frac{S_N}{\sqrt{n}}\right\} \right) d\tau, \qquad (5.3)$$

where $\sum'_{l,s}$ is the summation over all $l$-tuples $(s_1,...,s_l)$ with non-negative integer components such that $s_1 + ... + s_l = s$ ; $\sum''_l$ is the summation over all $l$-tuples $(j_1,...,j_l)$ such that $j_i \neq j_r$ for $i \neq r$ and $j_m = 1, 2, ..., N$ ; $m = 1, 2, ..., l$ .

Set $S_{l,N} = \sum_{i=1}^{l} (\xi_{j_i} - np_{j_i})$, $d_{l,N} = \sum_{i=1}^{l} p_{j_i}$ , $l = 1,...,s$ and write

$$\int_{-\pi\sqrt{n}}^{\pi\sqrt{n}} E\left( \left(h_{j_1}(\xi_{j_1})\right)^{s_1} \cdot \ldots \cdot \left(h_{j_l}(\xi_{j_l})\right)^{s_l} \exp\left\{i\tau \frac{S_N}{\sqrt{n}}\right\} \right) d\tau$$

$$= \int_{-\pi\sqrt{n}}^{\pi\sqrt{n}} E\left( \left(h_{j_1}(\xi_{j_1})\right)^{s_1} \cdot \ldots \cdot \left(h_{j_l}(\xi_{j_l})\right)^{s_l} \exp\left\{i\tau \frac{S_{l,N}}{\sqrt{n}}\right\} \right) E \exp\left\{i\tau \frac{S_N - S_{l,N}}{\sqrt{n}}\right\} d\tau$$

$$= \int_{-\pi\sqrt{n}}^{\pi\sqrt{n}} E\left( \left(h_{j_1}(\xi_{j_1})\right)^{s_1} \cdot \ldots \cdot \left(h_{j_l}(\xi_{j_l})\right)^{s_l} \exp\left\{i\tau \frac{S_{l,N}}{\sqrt{n}}\right\} \right) \left( E \exp\left\{i\tau \frac{S_N - S_{l,N}}{\sqrt{n}}\right\} - \exp\left\{-\frac{\tau^2}{2}(1 - d_{l,N})\right\} \right) d\tau$$

$$+ \int_{-\pi\sqrt{n}}^{\pi\sqrt{n}} \exp\left\{-\frac{\tau^2}{2}(1 - d_{l,N})\right\} E\left( \left(h_{j_1}(\xi_{j_1})\right)^{s_1} \cdot \ldots \cdot \left(h_{j_l}(\xi_{j_l})\right)^{s_l} \left( \exp\left\{i\tau \frac{S_{l,N}}{\sqrt{n}}\right\} - 1 \right) \right) d\tau$$

$$+ E\left[ \left(h_{j_1}(\xi_{j_1})\right)^{s_1} \cdot \ldots \cdot \left(h_{j_l}(\xi_{j_l})\right)^{s_l} \right] \int_{-\pi\sqrt{n}}^{\pi\sqrt{n}} \exp\left\{-\frac{\tau^2}{2}(1 - d_{l,N})\right\} d\tau \stackrel{def}{=} J_1 + J_2 + J_3. \qquad (5.4)$$

We have

$$E \exp\left\{\frac{i\tau(\xi_m - np_m)}{\sqrt{n}}\right\} = \exp\left\{np_m\left(e^{it/\sqrt{n}} - 1 - \frac{it}{\sqrt{n}}\right)\right\} = \exp\left\{-\frac{\tau^2}{2}p_m + \frac{\theta\tau^3}{6\sqrt{n}}p_m\right\}, \qquad (5.5)$$

$$\left| E \exp\left\{\frac{i\tau(\xi_m - np_m)}{\sqrt{n}}\right\} \right| = \exp\left\{-2np_m \sin^2 \frac{\tau}{2\sqrt{n}}\right\} \leq \exp\left\{-\frac{2p_m}{\pi^2}\tau^2\right\}, \qquad (5.6)$$

because $\sin^2 u/2 \geq u^2/\pi^2$, $|u| \leq \pi$, and

$$d_{lN} \leq k_n p_{\max} = o(1), \ l = 1,...,s, \qquad (5.7)$$

since (2.10).

In order to get an upper bound for the $|J_1|$ we first write the integral $J_1$ as a sum of two integrals, say $J'_1$ and $J''_1$, over intervals $|\tau| \leq \pi\sqrt{n}/2$ and $\pi\sqrt{n}/2 \leq |\tau| \leq \pi\sqrt{n}$, respectively. Next, noting that $E \exp\{i\tau(S_N - S_{l,N})/\sqrt{n}\} = \prod_m E \exp\{i\tau(\xi_m - np_m)/\sqrt{n}\}$ , where the product runs over $(1,...,N)$ except $(j_1,...,j_l)$, we use relations (5.5) and (5.6), respectively in $J'_1$ and $J''_2$. Then quite clear algebra gives

$$|J_1| \leq E|h_{j_1}(\xi_{j_1})|^{s_1} \cdot \ldots \cdot E|h_{j_l}(\xi_{j_l})|^{s_l} \left[ \int_{-\pi\sqrt{n}/2}^{\pi\sqrt{n}/2} \exp\left\{-\frac{\tau^2}{2}(1-d_{l,N})\right\} \left|\exp\left\{\frac{\theta\tau^3}{6\sqrt{n}}(1-d_{l,N})\right\} - 1\right| d\tau \right.$$

$$+ \int_{\frac{\pi\sqrt{n}}{2} \leq |\tau| \leq \pi\sqrt{n}} \left( \left| E \exp\left\{ i\tau \frac{S_N - S_{l,N}}{\sqrt{n}} \right\} \right| + \exp\left\{-\frac{\tau^2}{2}(1-d_{l,N})\right\} \right) d\tau \right]$$

$$\leq E|h_{j_1}(\xi_{j_1})|^{s_1} \cdot \ldots \cdot E|h_{j_l}(\xi_{j_l})|^{s_l} \left[ \int_{-\pi\sqrt{n}/2}^{\pi\sqrt{n}/2} \exp\left\{-\frac{6-\pi}{12}(1-d_{l,N})\tau^2\right\} \frac{|\tau|^3}{6\sqrt{n}}(1-d_{l,N}) d\tau \right.$$

$$+ C_5 \exp\left\{-\frac{n}{4}(1-d_{l,N})\right\} \right].$$

Hence,

$$J_1 = O\left(\frac{1}{\sqrt{n}}\right) E\left( (h_{j_1}(\xi_{j_1}))^{s_1} \cdot \ldots \cdot (h_{j_l}(\xi_{j_l}))^{s_l} \right). \tag{5.8}$$

since (5.7) and the functions $h_m(\cdot) \geq 0$. Next apply inequality $|e^{it} - 1 - it| \leq t^2/2$ to get

$$J_2 = \frac{C_6}{n} E\left( |h_{j_1}(\xi_{j_1})|^{s_1} \cdot \ldots \cdot |h_{j_l}(\xi_{j_l})|^{s_l} \left( \sum_{i=1}^{l} (\xi_{j_i} - np_{j_i}) \right)^2 \right)$$

$$\leq \frac{C_7}{n} l \sum_{i=1}^{l} E\left( |h_{j_i}(\xi_{j_i})|^{s_i} (\xi_{j_i} - np_{j_i})^2 \right) \prod_{m \neq i, m=1}^{l} E|h_{j_m}(\xi_{j_m})|^{s_m}$$

$$\leq C_8 E\left( (h_{j_1}(\xi_{j_1}))^{s_1} \cdot \ldots \cdot (h_{j_l}(\xi_{j_l}))^{s_l} \right) \frac{s}{n} \sum_{i=1}^{l} \left( s_i^{a_1} (np_{j_i})^{b_1} + s_i^{a_2} (np_{j_i})^{b_2} \right)$$

$$\leq C_9 E\left( (h_{j_1}(\xi_{j_1}))^{s_1} \cdot \ldots \cdot (h_{j_l}(\xi_{j_l}))^{s_l} \right) \frac{s}{n} \left( (np_{\max})^{b_1} s^{\max(1,a_1)} + (np_{\max})^{b_2} s^{\max(1,a_2)} \right)$$

$$= o(1) E\left( (h_{j_1}(\xi_{j_1}))^{s_1} \cdot \ldots \cdot (h_{j_l}(\xi_{j_l}))^{s_l} \right), \tag{5.9}$$

since the fact that $s_1^a + \ldots + s_l^a = s^{\max(1,a)}$ (because $s_1 + \ldots + s_l = s$), $s \leq k_n$ and conditions (2.9), (2.10). By a simple algebra we obtain

$$J_3 = \sqrt{2\pi} E\left( (h_{j_1}(\xi_{j_1}))^{s_1} \cdot \ldots \cdot (h_{j_l}(\xi_{j_l}))^{s_l} \right)(1+o(1)). \tag{5.10}$$

since (5.7). Now apply (5.8), (5.9) and (5.10) in the (5.4) to get

$$\int_{-\infty}^{\infty} E\left( (h_{j_1}(\xi_{j_1}))^{s_1} \cdot \ldots \cdot (h_{j_l}(\xi_{j_l}))^{s_l} \exp\left\{ i\tau \frac{S_N}{\sqrt{n}} \right\} \right) d\tau$$

$$= \sqrt{2\pi}(1+o(1)) E\left( (h_{j_1}(\xi_{j_1}))^{s_1} \cdot \ldots \cdot (h_{j_l}(\xi_{j_l}))^{s_l} \right). \tag{5.11}$$

Note that in (5.9) and (5.10), and hence in (5.11), the $o(1)$ is uniform in all $l$-tuples $(s_{j_1},...,s_{j_l})$ and $s \leq k_n$. For every integer $s \leq k_n$ and $k_n$ satisfying the condition (2.10) the relations (5.2), (5.3) and (5.11) imply $ER_N^s = ET_N^s(1+o(1))$. Proposition 5.1 follows from this equality and Assertion A1.4.

Let $\lfloor x \rfloor$ stands for the largest integer which less than or equal to $x$, $\xi \sim Poi(\lambda)$, $\lambda > 0$, $\mu_v(\lambda) = E(\xi - \lambda)^v$ and $D_\lambda$ denote differentiation w.r.t. $\lambda$.

**Lemma 5.2.** For any integer $v \geq 2$ one has

$$\mu_v(\lambda) = v! \sum_{l=1}^{\lfloor v/2 \rfloor} c_{l,v} \lambda^l, \tag{5.12}$$

where

$$0 < c_{l,v} < 1/l!, \quad l = 1, 2, ..., \lfloor v/2 \rfloor, \tag{5.13}$$

and

$$(v+1)c_{l,v+1} = lc_{l,v} + c_{l-1,v-1}, \tag{5.14}$$

here $l = 1, 2, ..., \lfloor (v+1)/2 \rfloor$ if $v$ is an even, and $l = 1, 2, ..., \lfloor (v+1)/2 \rfloor - 1$,

$(v+1)c_{\lfloor (v+1)/2 \rfloor, v+1} = c_{\lfloor (v-1)/2 \rfloor, v-1}$ if $v$ is an odd; here we put $c_{0,v-1} = 0$.

**Proof.** Apply the Bruno's formula for $Ee^{i\tau(\xi-\lambda)} = \exp\{\lambda(e^{i\tau} - 1 - i\tau)\}$ to get

$$\mu_v(\lambda) = v! \sum \lambda^{s_2+...+s_v} \prod_{m=2}^{v} \frac{1}{s_m!(m!)^{s_m}} = v! \sum_{l=1}^{\lfloor v/2 \rfloor} c_{l,v} \lambda^l,$$

where $\sum$ is the summation over all non-negative $s_2,...,s_v$ such that $2s_2 + ... + vs_v = v$ and $l = s_2 + ... + s_v$, $c_{l,v} = \sum \prod_{m=2}^{v} \frac{1}{s_m!(m!)^{s_m}}$, and hence $0 < c_{l,v} < 1/l!$. On the other hand by Riordan (1936, eq.(4.8))

$$\mu_{v+1}(\lambda) = v\lambda\mu_{v-1}(\lambda) + \lambda D_\lambda \mu_v(\lambda). \tag{5.15}$$

Just to keep notation simple we put $\kappa = \lfloor v/2 \rfloor$ and $\tilde{c}_{l,v} = v!c_{l,v}$. Note that $\tilde{c}_{1,v} = 1$ for any integer $v \geq 0$. Let $v$ is even, then $\lfloor (v+1)/2 \rfloor = \kappa$ and $\lfloor (v-1)/2 \rfloor = \lfloor (v-2)/2 \rfloor = \kappa - 1$. Therefore (5.12) and (5.15) imply

$$\sum_{l=1}^{\kappa} \tilde{c}_{l,v+1}\lambda^l = \sum_{l=1}^{\kappa-1} v\tilde{c}_{l,v-1}\lambda^{l+1} + \sum_{l=1}^{\kappa} l\tilde{c}_{l,v}\lambda^l = \sum_{l=2}^{\kappa}\left(v\tilde{c}_{l-1,v-1} + l\tilde{c}_{l,v}\right)\lambda^l + \lambda.$$

For the even $v$ the property (5.14) follows. Let now $v$ is odd, then $\lfloor (v-1)/2 \rfloor = \kappa$ and $\lfloor (v+1)/2 \rfloor = \kappa + 1$. From (5.12) and (5.15) obtain

$$\sum_{l=1}^{\kappa+1} \tilde{c}_{l,v+1}\lambda^l = \sum_{l=1}^{\kappa} v\tilde{c}_{l,v-1}\lambda^{l+1} + \sum_{l=1}^{\kappa} l\tilde{c}_{l,v}\lambda^l = \sum_{l=2}^{\kappa+1} v\tilde{c}_{l-1,v-1}\lambda^l + \sum_{l=1}^{\kappa} l\tilde{c}_{l,v}\lambda^l$$

$$= \sum_{l=2}^{\kappa} \left( v\tilde{c}_{l-1,v-1} + l\tilde{c}_{l,v} \right) \lambda^l + v\tilde{c}_{\kappa,v-1} \lambda^{\kappa+1} + \lambda.$$

This equality proves (5.14) for the odd $v$. Lemma 5.2 is proved completely.

Using (5.12) we obtain $\mu_v(\lambda) < \lambda D_\lambda \mu_v(\lambda) \leq \lfloor v/2 \rfloor \mu_v(\lambda)$. This together with (5.15) imply

$$\mu_v(\lambda) < \mu_{v+1}(\lambda) \leq \left( v\lambda + \lfloor v/2 \rfloor \right) \mu_v(\lambda). \tag{5.16}$$

Use formula (5.15) for $\mu_{v+2}(\lambda)$, next again apply (5.15) to the derivative of $\mu_{v+1}(\lambda)$, after this use the inequalities $\lambda D_\lambda \mu_{v-1}(\lambda) \leq \mu_v(\lambda)$ (because $l\tilde{c}_{l,v-1} \leq \tilde{c}_{l,v}$, see (5.14)), $\lambda D_\lambda \mu_V(\lambda) \leq \lfloor v/2 \rfloor \mu_V(\lambda)$ and $\lambda^2 D_\lambda^2 \mu_v(\lambda) \leq \lfloor v/2 \rfloor^2 \mu_v(\lambda)$. Then one can observe that for any integer $v \geq 0$

$$(v+1)\lambda \mu_v(\lambda) < \mu_{v+2}(\lambda) \leq 2v(\lambda + v)\mu_v(\lambda). \tag{5.17}$$

**Proof of Theorem 2.2.** Due to Proposition 5.1 and condition (ii) we have for large enough $N$

$$\left| \mathcal{C}_s(\sigma_N^{-1} R_N(\eta)) \right| \leq \sigma_N^{-s} \sum_{m=1}^{N} \left| \mathcal{C}_s(h_m(\xi_m)) \right| \leq (s!)^{1+v} V_n^{s-2} \tilde{\sigma}_N^2 \sigma_N^{-s}$$

$$\leq (s!)^{1+v} (\sigma_N V_n^{-1})^{-(s-2)} (\sigma_N^2 / \tilde{\sigma}_N^2)^{-1} \leq (s!)^{1+v} (V_n^{-1} \sigma_N^3 / \tilde{\sigma}_N^2))^{1/(s-2)})^{-(s-2)} = (s!)^{1+v} W_N^{-(s-2)}.$$

for all $s$: $3 \leq s \leq k_n$, see (2.10), because $\tilde{\sigma}_N^2 \geq \sigma_N^2$. Theorem 2.2 follows now from Assertion 2 (see Appendix) with $\Delta = W_N$ and $\vartheta = k_n$.

**Proof of Theorem 3.1.** From (3.3), (3.4) for $d \neq 0$, and for $d = 0$ from (3.6), (3.7), and the fact that $h_{0,m}(\xi_m) = 2np_m \left(1 + \tilde{\xi}_m\right) \ln\left(1 + \tilde{\xi}_m\right)$, where $\tilde{\xi}_m = (\xi_m - np_m)/np_m$, using inequalities $x \leq (1+x)\ln(1+x) \leq x + x^2/2$ one can observe that under the condition (3.9) for large enough $N$

$$\tilde{\sigma}_N^2(d) \sim c_1 N \text{ and } \sigma_N^2(d) \sim c_2 N. \tag{5.18}$$

Case (i) follows from Theorem 2.1, (5.18) and that (see (3.1))

$$E \exp\left\{ H \left| h_{d,m}(\xi_m) \right| \right\} \leq E \exp\left\{ 2H(np_m)^{-d} \xi_m \right\} \leq \exp\left\{ np_m \left( \exp\left\{ 2H(np_m)^{-d} \right\} - 1 \right) \right\} \leq C(d, H, c_3, c_4),$$

since condition (3.9).

**Proof** of the cases (ii) consist in verifying of conditions of Theorem 2.2. It is easy to verify that $E\xi^{s+1} = \lambda(E\xi^s + D_\lambda E\xi^s)$, $\xi \sim Poi(\lambda)$, for arbitrary $s \geq 0$. Using this and raw inequality $\lambda^i D_\lambda^i E\xi^s < s^i E\xi^s$ it is not hard to see that $E\xi^{s+2} \leq 4\max(s^2, \lambda^2) E\xi^s$ for $s \geq 0$. Hence, condition (2.9) is fulfilled with $a_1 = 2, b_1 = 0$, for $s > \lambda$, and $a_2 = 0, b_2 = 2$, for $s \leq \lambda$, that is

$$K(a_1, b_1) = n^{1/3} = c_1 N^{1/3} \text{ and } K(a_2, b_2) = (n^{-1} p_{\max}^{-2})^{1/2} = c_2 N^{1/2}. \tag{5.19}$$

Further, by Bruno's formula for derivatives of composite function (Petrov (1993, Eq. (1.6)) and Eq (5.2) of Sachkov (1996), we have for any integer $k \geq 1$

$$E\xi_m^k = k! \sum_{l=1}^{k} c_{l,k}(np_m)^l, \ 0 < c_{l,k} \leq 1. \tag{5.20}$$

Next, using inequalities $\sqrt{2\pi k}(k/e)^k < k! < e\sqrt{2\pi k}(k/e)^k$ by simple algebra we obtain: for $d > 0$

$$(\lfloor s(d+1) \rfloor + 1)! \leq (s!)^{1+d}[C(d)]^{s-2}, \ C(d) = e(2+d)^{3(1+d)}(3(1+d)+1)/3^{d/2}. \tag{5.21}$$

We have $Varh_{d,m}(\xi_m) \sim c$, since condition (3.9). By this fact and (5.20), (5.21) we obtain

$$E\left(h_{d,m}(\xi_m) - Eh_{d,m}(\xi_m)\right)^s \leq 2^s(np_m)^{-ds} E\xi_m^{s(d+1)} \leq 2^s(np_m)^{-ds} E\xi_m^{\lfloor s(d+1) \rfloor + 1}$$

$$\leq (s!)^{1+d}[C_0(d,\lambda,c_3,c_4)]^{s-2} Varh_{d,m}(\xi_m).$$

Hence by Assertion 1

$$\left|\mathcal{C}_s(h_{d,m}(\xi_m))\right| \leq (s!)^{1+d}[4C_0(d,c_3,c_4)]^{s-2} Varh_{d,m}(\xi_m). \tag{5.22}$$

Thus, condition (2.11) is fulfilled with $v = d > 0$, $V_n = 4C_0(d,c_3,c_4)$ and hence $W_n = \sqrt{N}$, since (5.18). Due to these, (5.19) and (2.12) we have $x_n = o(N^{1/6})$, therefore $M_{d,n}(x_n) = o(1)$, since (2.15). Case (ii) follows. The proof of Theorem 3.1 is complete.

**Proof of Theorem 3.2.** The PDS can be rewritten as

$$CR_N(d) = \sum_{l=1}^{N} np_j \tilde{h}_{d,j}(\eta_j),$$

where

$$\tilde{h}_{d,j}(x) = \frac{2}{d(d+1)}\left[(x/np_l)^{d+1} - (d+1)(x/np_l) + d\right]$$

$$= \frac{2}{d(d+1)}\left[\left(1 + \frac{x - np_j}{np_l}\right)^{d+1} - (d+1)\left(1 + \frac{x - np_j}{np_l}\right) + d\right], \ d \neq 0,$$

$$\tilde{h}_{0,j}(x) = 2\left[(x/np_l)\log(x/np_l) - (x - np_j)/np_j\right]$$

$$= 2\left[\left(1 + \frac{x - np_j}{np_j}\right)\log\left(1 + \frac{x - np_j}{np_j}\right) - \frac{x - np_j}{np_j}\right]. \tag{5.23}$$

Set $\tilde{\xi}_m = (\xi_m - np_m)/np_m$. Recall that $np_{\min} \to \infty$, therefore the r.v. $\sqrt{np_m}\tilde{\xi}_m$ has asymptotically normal distribution. Hence

$$\tilde{\xi}_m = O_p\left((np_m)^{-1/2}\right). \tag{5.24}$$

Now, we apply Taylor expansion formula in (5.23) to get

$$\tilde{h}_{d,m}(\xi_m) = \tilde{\xi}_m^2 + O_p\left(\tilde{\xi}_m^3\right) = \tilde{\xi}_m^2 + O_p\left((np_m)^{-3/2}\right), \ d > -1. \tag{5.25}$$

By virtue of this fact and that $E\tilde{\xi}_m^2 = (np_m)^{-1}$, $EO_p\left((np_m)^{-3/2}\right) = O\left((np_m)^{-3/2}\right)$ we obtain under condition (3.12) for all $d > -1$: $A_N(d) = N + o(N)$ and

$$\sigma_N^2(d) = \sum_{m=1}^{N} \frac{1}{np_m} + (2 - \lambda_n^{-1})N + \sum_{m=1}^{N} O\left((np_m)^{-5/2}\right) = 2N(1 + o(1)), \quad (5.26)$$

since $\lambda_n \to \infty$ as $np_{\min} \to \infty$, because $1 = p_1 + ... + p_N \geq Np_{\min}$ and hence $np_{\min} \leq \lambda_n$. Further proof consists in verifying the conditions of Theorem 2.2. Note that $\tilde{h}_{d,m}(x) \geq 0$ for all $d \in (-\infty, \infty)$. Next, by (5.24) and (5.25) we have

$$E\tilde{h}_{d,m}^s(\xi_m) = E\tilde{\xi}_m^{2s}\left(1 + O_p\left((np_m)^{-1/2}\right)\right) = E\tilde{\xi}_m^{2s} + O\left((np_m)^{-(2s+1)/2}\right), \quad (5.27)$$

since $EO_p\left((np_m)^{-(2s+1)/2}\right) = O\left((np_m)^{-(2s+1)/2}\right)$ and by Lemma 5.2 $E\tilde{\xi}_m^{2s} > c(s)(np_m)^{-s}$, where $1 \leq c(s) < s$. Similarly, by (5.24), (5.25) and (5.17) we get

$$E\tilde{h}_{d,m}^s(\xi_m)(\xi_m - np_m)^2 = E\tilde{\xi}_m^{2s}(\xi_m - np_m)^2 + O\left((np_m)^{-(2s-1)/2}\right)$$

$$\leq 4s(np_m + 2s)E\tilde{\xi}_m^{2s} + O\left((np_m)^{-(2s-1)/2}\right). \quad (5.28)$$

So condition (2.9) is fulfilled with $b = 1$ and $a = 1$ if $np_m \geq 2s$, and with $b = 0$ and $a = 2$ if $np_m < 2s$. That is

$$K(a_1, b_1) = p_{\max}^{-1/2} \text{ and } K(a_2, b_2) = n^{1/3}. \quad (5.29)$$

Next, by (5.25) we obtain

$$E(np_m \tilde{h}_{d,m}(\xi_m)) = 1 + O\left((np_m)^{-1/2}\right), \quad Var(np_m \tilde{h}_{d,m}(\xi_m)) = 2 + O\left((np_m)^{-1/2}\right) \quad (5.30)$$

Use (5.25), (5.30) and Lemma 3.2 (relations (5.12) and (5.13)) to get

$$\left|E\left(np_m\left(\tilde{h}_{d,m}(\xi_m) - E\tilde{h}_{d,m}(\xi_m)\right)\right)^s\right| \leq 2^s (np_m)^s E\tilde{h}_{d,m}^s(\xi_m)$$

$$\leq 2^{2s}\left((np_m)^{-s}E(\xi_m - np_m)^{2s} + O\left((np_m)^{-s/2}\right)\right)$$

$$< (2s)! 2^{2s} e \leq (s!)^2 2^{4s} e^2 / \pi\sqrt{s} \leq (s!)^2 2^{8(s-2)} Var(np_m \tilde{h}_{d,m}(\xi_m)).$$

since $(2s)! \leq (s!)^2 2^{2s} e / \pi\sqrt{2s}$, due to inequalities $\sqrt{2\pi m}\, m^m e^{-m} \leq m! \leq e\sqrt{m}\, m^m e^{-m}$. Therefore, from Assertion1 it follows that

$$\left|\mathcal{C}_s\left(np_m \tilde{h}_{d,m}(\xi_m)\right)\right| \leq (s!)^2 2^{9(s-2)} Var(np_m \tilde{h}_{d,m}(\xi_m)).$$

Hence, condition (2.11) is fulfilled with $V_n = 2^9$, $v = 1$ and hence $M_{v,n}(x) = 0$. Yet $W_n = \sqrt{N}$, because (5.26) and $\tilde{\sigma}_N^2 = 2N(1 + o(1))$ by (5.30). Theorem 3.2 follows.

**Proof of Theorems 3.3 and 3.4** will be concluded by applying Theorem 2.2. Remind $np_{\max} \to 0$. Due to this fact for arbitrary well defined function $\varphi(x)$ we obtain

20$$E\varphi(\xi_m) = \sum_{l=0}^{\infty}(-1)^l \frac{(np_m)^l}{l!}\sum_{k=0}^{\infty}\varphi(k)\frac{(np_m)^k}{k!}$$

$$= \sum_{v=0}^{\infty}\frac{(np_m)^v}{v!}\sum_{l=0}^{v}\binom{v}{l}(-1)^l\varphi(v-l) = \sum_{v=0}^{\infty}\frac{(np_m)^v}{v!}\Delta^v\varphi(0), \quad (5.31)$$

where $\Delta^0\varphi(x) = \varphi(x)$, $\Delta\varphi(x) = \varphi(x+1) - \varphi(x)$. The Eq.s (3.14) - (3.19) can be derived from (5.31) by not hard algebra. For $h_{d,l}(x)$ defined in (2.1) from (5.31) we obtain: for $d \ne 0$

$$E(h_{d,m}(\xi_m))^s = \sum_{v=1}^{\infty}\frac{(np_m)^v}{v!}\Delta^v h_{d,m}^s(0),$$

and

$$E(h_{d,m}(\xi_m))^s(\xi_m - np_m)^2 = (np_m)^{-sd}\sum_{v=1}^{\infty}\frac{(np_m)^v}{v!}\sum_{l=0}^{v}\binom{v}{l}(-1)^l(v-l)^{s(1+d)}(v-l-np_m)^2$$

$$= \sum_{v=1}^{\infty}\frac{(np_m)^v}{v!}v^2(1+o(1))\Delta^v h_{d,m}^s(0) = O(Eh_{d,m}^s(\xi_m)).$$

Here the $O(.)$ is uniformly in $s$. Similar result for the case $d = 0$ follows from (5.31) with $\varphi(x) = x\log(x/np_m)$. Thus, the condition (3.2) of Lemma 3.1 is fulfilled with $a = 0$ and $b = 0$, i.e. $K_n(a,b) = \sqrt{n}$.

The next step of the proofs consist in use of Assertion 2, see Appendix. To end this we have to get an appropriate inequalities for the cumulants of $h_{d,m}(\xi_m)$. Since $np_{\max} \to 0$ we obtain

$$E\xi_m^{s(d+1)} \le E\xi_m^{\lfloor s(d+1)\rfloor+1} \le (\lfloor s(d^*+1)\rfloor+1)!(np_m)(1+o(1)), \quad (5.32)$$

where $d^* = \max(0,d)$. Alike (5.21) we obtain

$$(\lfloor s(d^*+1)\rfloor+1)! \le (s!)^{1+d^*}[C(d^*)]^{s-2}. \quad (5.33)$$

We have

$$Var\xi_m^{1+d} = np_m(1+O(np_m)). \quad (5.34)$$

*Proof of Theorem 3.3.* Assume $(p_1,...,p_N) \ne (N^{-1},...,N^{-1})$. Let $d \ne 0$. By (5.32), (5.33) and (5.34) we get

$$\left|E\left(h_{d,m}(\xi_m) - Eh_{d,m}(\xi_m)\right)^s\right| \le 2^s(np_m)^{-ds}E\xi_m^{s(d+1)} \le 2^s(np_m)^{-ds}E\xi_m^{\lfloor s(d^*+1)\rfloor+1}$$

$$\le 2(\lfloor s(d^*+1)\rfloor+1)!(np_m)^{1-ds} \le (s!)^{1+d^*}[C(d)]^{s-2}(np_m)^{1-2d}(np_m)^{-d(s-2)}$$

$\le (s!)^{1+d^*}[C(d)(np_m)^{-d}]^{s-2}(np_m)^{1-2d} \le (s!)^{1+d^*}[2C(d)(np_m)^{-d}]^{s-2}Varh_{d,m}(\xi_m)$. Hence, by Assertion 1,

$$\left|\mathcal{C}_s(h_{d,m}(\xi_m))\right| \le (s!)^{1+d^*}[4C(d)(np_m)^{-d}]^{s-2}Varh_{d,m}(\xi_m), d \ne 0. \quad (5.35)$$

because $\sqrt{Varh_{d,m}(\xi_m)} = \sqrt{np_m}(np_m)^{-d} = o((np_m)^{-d})$. Hence, condition (2.11) is fulfilled with

$$V_n = 4C(d)(np_m)^{-d}, v = d^*, W_N = \tilde{\sigma}_N^{-2}(d)\sigma_N^3(d)(np_{\min})^d. \qquad (5.36)$$

Let now $d = 0$, that is $h_{0,m}(x) = 2x\log x - 2x\log(np_m)$. Since $np_{\max} \to 0$, we find that

$$Varh_{0m}(\xi_m) = 4np_m \ln^2 np_m + 8(np_m)^2[\ln^2 2 - 2\ln 2\ln np_m](1+o(1)) = 4np_m \ln^2 np_m (1+o(1)).$$

Next, we have

$$\left|E(h_{0,m}(\xi_m) - E(h_{0,m}(\xi_m)))^s\right| \le 2^{2s} E(\xi_m \log(\xi_m/np_m))^s$$

$$\le 2^{3s-1} E(\xi_m^s \ln^s \xi_m + \xi_m^s(-\ln np_m)^s)) \le 2^{3s-1}[E\xi_m^{s+1} + (-\ln np_m)^s E\xi_m^s]$$

$$\le 2^{3s}[(s+1)!(np_m) + s!(-\ln np_m)^s(np_m)] \le s!(-8\ln np_m)^s(np_m)$$

$$\le s![-2^7 \ln np_m]^{s-2} 4(np_m)\ln^2 np_m \le s![-2^7 \ln np_{\min}]^{s-2} Varh_{0m}(\xi_m).$$

From this and Assertion 1 it follows that

$$\left|\mathcal{C}_s(h_{0,m}(\xi_m))\right| \le s![2|\ln np_{\min}|]^{s-2} Varh_{0m}(\xi_m).$$

Hence, condition (2.11) is fulfilled with $V_n = 2|\ln np_{\min}|, v = 0, W_N = \sigma_N^3(0)/\tilde{\sigma}_N^{-2}(0)|\ln np_{\min}|$. By this, (5.36) and that $K_n(a,b) = \sqrt{n}$ Theorem 3.3 follows from Theorem 2.2.

*Proof of Theorem 3.4.* Now $h_{d,m}(x) = h_d(x) = x^{1+d}$, $d \ne 0$, $h_0(x) = 2x\ln x$, and $\xi_m \sim Poi(\lambda_n)$. Therefore, alike (5.35) we obtain

$$\left|\mathcal{C}_s(h_d(\xi_m))\right| \le (s!)^{1+d^*}[2C(d)]^{s-2}\lambda_n \le (s!)^{1+d^*}[4C(d)]^{s-2} Varh_d(\xi_m), \qquad (5.37)$$

since $Varh_d(\xi_m) = \lambda_n(1+o(1))$, see (5.34). Hence in this case by (3.20)

$$v = d^*, W_N = \sqrt{n}\lambda_n^{3/2}(2^d - 1)^3/\sqrt{2}C(d). \qquad (5.38)$$

Further, we have

$$E(h_0(\xi_m))^s = E(\xi_m \log \xi_m)^s \le E\xi_m^{s+1} \le (s+1)!\lambda_n(1+o(1))$$

$$\le s![((s+1)/\lambda_n 8\ln^2 2)^{1/(s-2)}]^{s-2}\lambda_n^2 8\ln^2 2 \le s![(\lambda_n 2\ln^2 2)^{-1}]^{s-2} Varh_0(\xi_m),$$

since $Varh_0(\xi_m) = 8\ln^2 2 \lambda_n^2$. Hence by Assertion 1

$$\left|\mathcal{C}_s(h_0(\xi_m))\right| \le s![8(\lambda_n \ln^2 2)^{-1}]^{s-2} Varh_0(\xi_m). \qquad (5.39)$$

Hence, $v = 0$, $W_N = 16\sqrt{2}\ln^{5/2} 2\sqrt{n}\lambda_n^{3/2}$, by (3.19). Theorem 3.4 follows from Theorem 2.2, since (5.37), (5.38), (5.39) and that $K_n(a,b) = \sqrt{n}$.

**Proof of Corollaries 3.1, 3.2 and 3.3** follows from Theorems 3.1, 3.2 and 3.4, respectively, and that $\log(1-\Phi(x)) \sim -x^2/2$, as $x \to \infty$, and by (2.15).

**Proof of Theorems 4.1, 4.2 and 4.3** follows in immediate manner from Theorem 3.1 (ii) case $d = 1$, Theorem 3.2, 3.3 and 3.4, respectively.

**Proof of Theorem 4.4**. For the chi-square statistic $h_m(x) = (x - np_m)^2 / np_m$, hence $A_N = N$, $\sigma_N^2$ and $\tilde{\sigma}_N^2$ as in (4.1). Set $\tilde{\xi}_m = (\xi_m - np_m)/\sqrt{np_m}$. Then $h_m(\xi_m) = \tilde{\xi}_m^2$ and $E\tilde{\xi}_m^{2s}(\xi_m - np_m)^2 \leq 4s(np_m + 2s)E\tilde{\xi}_m^{2s}$, since (5.17), hence the condition (2.9) is fulfilled with $b = 1$ and $a = 1$ if $np_m \geq 2s$, and with $b = 0$ and $a = 2$ if $np_m < 2s$. We would remind the notation $\nabla_n = \max\{1,(np_{\min})^{-1}\}$. Set $\varsigma_m = (\tilde{\xi}_m^2 - E\tilde{\xi}_m^2)/\sqrt{Var\tilde{\xi}_m^2}$. Using Lemma 5.2 we obtain

$$\left|E\varsigma_m^s\right| \leq 2^s (2s)!(Var\tilde{\xi}_m^2)^{-s/2} \sum_{l=0}^{s-1}(np_m)^{-l} \leq (s!)^2 \left(2^7 \nabla_n / \sqrt{Var\tilde{\xi}_m^2}\right)^{s-2} \text{ for all } s \geq 3,$$

since $(2s)! \leq (s!)^2 2^{2s} e / \pi\sqrt{2s}$, $\sum_{l=0}^{s-1}(np_m)^{-l} \leq s\nabla_n^{s-2}$ and $Var\tilde{\xi}_m^2 = (np_m)^{-1} + 2$. Hence by Assertion 1

$$\left|\mathcal{C}_s(\varsigma_m)\right| \leq (s!)^2 \left(2^8 \nabla_n / \sqrt{Var\tilde{\xi}_m^2}\right)^{s-2}. \text{ That is } \left|\mathcal{C}_s(\tilde{\xi}_m^2)\right| \leq (s!)^2 \left(2^8 \nabla_n\right)^{s-2} Var\tilde{\xi}_m^2, \ s \geq 3. \text{ Thus condition}$$

(ii) of Theorem 2.2 is fulfilled with $v = 1$, $V_n = 2^8 \nabla_n$ and $3 \leq s \leq k_n = o\left(\min(n^{1/3}, p_{\max}^{-1/2})\right)$, because $K_n(1,1) = p_{\max}^{-1/2}$ and $K_n(2,0) = n^{1/3}$. Theorem 4.2 follows from Theorem 2.2.

**Proof of Corollary 4.1**. We have $\tilde{\sigma}_N^2 \leq N(2 + c_0^{-1})$, $\nabla_n \leq \max(1, c_0^{-1})$ and $n \geq c_0 N$, since $np_{\min} \geq c_0$. Yet $\sigma_N^2 \geq 2N$, see Remark 4.1. Corollary 2.1 follows.

**Proof of Corollary 4.2** is straightforward, since in this case $\nabla_n = \max(1, \lambda_n^{-1})$, $\sigma_N^2 = N\left(2 + d^2 \delta^2(n)\lambda_n^{-1}(1 + o(1))\right)$ and $\tilde{\sigma}_N^2 = 2N\left(1 + (2\lambda_n)^{-1}(1 + o(1))\right) \leq 2N\nabla_n$.

**Proof of Theorem 4.5 and 4.6** follows straightforwardly from Theorems 3.1 and 3.2 by putting $d = 0$ and $d = -1/2$, respectively.

**Proof of Theorem 4.7** follows from Theorem 2.2. Indeed. In this case $h_m(x) = \mathrm{I}\{x = r\}$ and it is easy to see that the condition (i) of Theorem 2.2 is fulfilled with $\mathcal{N} = \{j : np_j > r\}$

$$a_1 = 0, \ b_1 = 2 \text{ and } a_2 = 0, b_2 = 0. \tag{5.40}$$

Further. It is obvious that $Varh_m(\xi_m) = \pi_r(np_m)(1 - \pi_r(np_m))$ and $Eh_m^s(\xi_m) = \pi_r(np_m)$, $s \geq 1$. We have $\left|E(h_m(\xi_m) - \pi_r(np_m))^s\right| \leq 2^s \pi_r(np_m) < s![8(1 - \pi_r(np_m))^{-1}]^{s-2} Varh_m(\xi_m) \leq s![8c(r)]^{s-2} Varh_m(\xi_m)$, $s \geq 3$ and a constant $c(r) > 0$. Therefore, by Assertion 1

$$\left|\mathcal{C}_s(h_m(\xi_m))\right| \leq s!(16c(r))^{s-2} Varh_m(\xi_m).$$

Hence condition (2.11) is fulfilled with $v = 0$, $V_n = 8c(r)$. Thus $W_n = \sigma_{rN}^3 / \tilde{\sigma}_{rN}^2$, $K_n(a_2, b_2) = n^{1/2}$ and $K_n(a_1, b_1) = (np_{\max}^2)^{-1/2}$, since (5.40). Theorem 4.7 follows.

**Proof of Corollary 4.3.** The case (i) follows from Theorem 2.1, because $E\exp\{H\mathrm{I}\{\xi_m=r\}\}\le e$ and (4.7). In the cases (ii) by (4.8) we have $\sigma_{rN}^3/\tilde{\sigma}_{rN}^2 = n^{r/2}/N^{(r-1)/2}$, $r\ge 2$, and $\sigma_{rN}^3/\tilde{\sigma}_{rN}^2 = n^{3r/2}/N^{(3r-1)/2}$, $r=0,1$. Therefore, it is seen that $\Upsilon_{r,n} = n^{1/4}$. So the case (ii) follows from Theorem 4.7. Case (iii) follows from Theorem 4.7 and (4.9) straightforwardly.

**Proof of Theorem 4.8.** Remind that for the statistic $\Phi_N$ the kernel function is $h_m(x) = \mathrm{I}\{v_m > x\}$. If $\lambda_n \to \lambda \in (0,\infty)$, then it is not hard to verify that all conditions of Theorem 2.1 are fulfilled, and hence the case (i) follows.

The cases (ii) and (iii) we will derive from Theorem 2.2. Observe that $Eh_m^s(\xi_m) = \tau(\lambda_n)$, $s\ge 1$. Using this we find that $Eh_m^s(\xi_m)(\xi_m - \lambda_n)^2 \le \lambda_n(3\lambda_n + 2)Eh_m^s(\xi_m)$. Hence the condition (i) of Theorem 2.2 is fulfilled with

$$a = 0, b = 0 \text{ if } \lambda_n \le 1, \tag{5.41}$$

and

$$a = 0, \ b = 2 \text{ if } \lambda_n > 1. \tag{5.42}$$

Next, alike to above proof of Theorem 4.7 it can be shown that

$$\left|\mathcal{C}_s(h_m(\xi_m))\right| \le s!(16(1-\tau(\lambda_n))^{-1})^{s-2} Varh_m(\xi_m), s\ge 3,$$

since $Varh_m(\xi_m) = \tau(\lambda_n)(1-\tau(\lambda_n))$. That is condition (2.11) is fulfilled with $v=0$, $V_n = 16(1-\tau(\lambda_n))^{-1}$. Hence

$$\Upsilon_n = \sqrt{N}\left(\tau(\lambda_n)(1-\tau(\lambda_n)) - \lambda_n(\tau'(\lambda_n))^2\right)^{3/2}/\tau(\lambda_n). \tag{5.43}$$

Let $\lambda_n \to 0$. Then from (4.10) and (4.11) by simple algebra we obtain $A_N = N(1-\beta_0)(1+o(1))$, $\sigma_N^2 = N\beta_0(1-\beta_0)(1+O(\lambda_m^\vartheta))$, $\tau(\lambda_n)(1-\tau(\lambda_n)) - \lambda_n(\tau'(\lambda_n))^2 = \beta_0(1-\beta_0)(1+O(\lambda_m^\vartheta))$. The case (ii) follows from Theorem 2.2, (5.41) and (5.43).

Case (iii). In this case from (4.6) we have $\sigma_N^2 \sim \tilde{\sigma}_N^2 = N\tau(\lambda_n)(1+o(1))$. If (4.12) is fulfilled and $\lambda_n \le \delta\mu^{-1}\ln N + \chi\ln\ln N$ then $\tau(\lambda_n) \le \rho N^{-\delta}$. Also, $\sigma_N^2 = Ng\lambda_n^{-\varepsilon}(1+o(1))$ if (4.13) is fulfilled. Since these facts and (5.42) and (5.43) case (iii) follows from Theorem 2.2.

**Appendix.**

Below $\varsigma$ is a r.v. with $E\varsigma = 0$, $Var\varsigma = \sigma^2 > 0$, $\mathcal{C}_k(\varsigma)$ and $\alpha_k(\varsigma)$ cumulant and moment of $k$ th order, respectively of the r.v. $\varsigma$. The following conditions and assertions are taken from the book by Saulis and Statulevicius (1991), SS(1991) for short.

**Bernstein condition** $(B_\nu)$: there exist constants $\nu \geq 0$ and $B > 0$ such that

$$|\alpha_k(\varsigma)| \leq (k!)^{\nu+1} B^{k-2} \sigma^2, \text{ for all } k = 3, 4, \ldots.$$

**Statulevicius condition** $(S_\nu)$: there exist constants $\nu \geq 0$ and $\Delta > 0$ such that

$$|\mathcal{C}_k(\varsigma)| \leq (k!)^{\nu+1} \Delta^{-(k-2)} \sigma^2, \text{ for all } k = 3, 4, \ldots.$$

**Assertion 1** (SS (1991, Lemma 3.1)). If $\varsigma$ satisfy condition $(B_\nu)$ then it also satisfy condition $(S_\nu)$ with $\Delta = (2\max(1, B))^{-1}$.

**Assertion 2.** Assume $\sigma^2 = 1$ and there exist $\nu \geq 0$, $\Delta > 0$ and $\kappa \geq 3$ such that

$$|\mathcal{C}_k(\varsigma)| \leq (k!)^{\nu+1} \Delta^{-(k-2)}, \text{ for all } k = 3, 4, \ldots, \kappa. \tag{A.1}$$

Then there exist constants $c_i > 0$ such that for all $x: 0 \leq x \leq c_1 \min(\sqrt{\kappa}, \Delta^{1/(1+2\nu)})$ it holds

$$P\{\xi > x\} = (1 - \Phi(x))\exp\{L_\nu(x)\}\left(1 + \theta_3 c_2 \frac{x+1}{\min(\sqrt{\kappa}, \Delta^{1/(1+2\nu)})}\right),$$

$$P\{\xi \leq -x\} = \Phi(-x)\exp\{L_\nu(-x)\}\left(1 + \theta_3 c_3 \frac{x+1}{\min(\sqrt{\kappa}, \Delta^{1/(1+2\nu)})}\right),$$

where $L_\nu(x) = x^3(l_0 + l_2 x + \ldots + l_9 x^9)$, $\vartheta = \vartheta(\nu): \vartheta(0) = \infty$, $\vartheta(\nu) < \nu^{-1} - 1$ if $\nu \in (0, 1)$, and $L_\nu(x) = 0$ if $\nu \geq 1$. Also $|L_\nu(x)| \leq 5|x|^3 / 4\Delta^{1/(1+2\nu)}$.

**Proof** based on a revision of the proofs of Lemma 2.1 and Lemma 2.2 from SS (1991). We state their Lemma 2.1 first.

**Lemma 2.1** of SS (1991). Let there exist $\Delta > 0$ and an integer $\upsilon \geq 3$ such that

$$|\mathcal{C}_k(\xi)| \leq (k-2)! \Delta^{-(k-2)}, \text{ for all } k = 3, 4, \ldots, \upsilon.$$

where an even $\upsilon \leq 2\Delta^2$, is fulfilled. Then for all $x: 0 \leq x \leq c_0 \sqrt{\upsilon}$, where $c_0 < (3\sqrt{e})^{-1}$, it holds

$$P\{\xi > x\} = (1 - \Phi(x))\exp\{L(x)\}\left(1 + \theta_1 f_1(x, \upsilon)(x+1)\upsilon^{-1/2}\right),$$

$$P\{\xi \leq -x\} = \Phi(-x)\exp\{L(-x)\}\left(1 + \theta_2 f_1(x, \upsilon)(x+1)\upsilon^{-1/2}\right),$$

where $f_1(x, \upsilon) = \left(1 - 3x\sqrt{e/\upsilon}\right)^{-1}\left(117 + 96\upsilon \exp\{-\upsilon^{1/4} 2^{-1}\left(1 - 3x\sqrt{e/\upsilon}\right)\}\right)$ and $L(x) = l_0 x^3 + l_1 x^4 + \ldots$

is a power series such that $|L(x)| \leq 5|x|^3 / 4\Delta$ and is converging within the disc $|x| < \Delta(3\sqrt{e}/2)^{-1}$.

As in inequality (2.61) of SS (1991) we have $(k!)^{1+\nu} \leq (k-2)!(6\kappa^\nu)^{k-2}$, for $k = 3, \ldots, \kappa$ and $\kappa \geq 6$. Therefore from condition (A.1) it follows

$$|\mathcal{C}_k(\xi)| \leq (k-2)!\Delta_s^{-(k-2)}, \text{ for all } k=3,...,\kappa, \text{ and } \Delta_\kappa = \Delta/6\kappa^\nu.$$

Let $\kappa < 2\Delta_\kappa^2$. If $\kappa$ is an even, then from Lemma 2.1 of SS(1991) gives

$$P\{\xi > x\} = (1-\Phi(x))\exp\{L(x)\}\left(1+\theta_1 c_2(x+1)\kappa^{-1/2}\right), \text{ for } 0 \leq x < c_1\sqrt{\kappa}, \qquad (A.2)$$

where $c_1 = (3\sqrt{e})^{-1}$ and $c_2 = f_1(x,\kappa)$. If $\kappa$ is an odd, then we use Lemma 2.1 putting $\upsilon = \kappa - 1$ to get (A.2) with $c_1 = \sqrt{5}/3\sqrt{6e}$ and $c_2 = \sqrt{6/5}\,f_1(x,\kappa-1)$.

Let, now, $\kappa \geq 2\Delta_\kappa^2$, that is $\kappa \geq \tilde{\Delta}_\nu := \left(\Delta^2/18\right)^{1/(1+2\nu)}$. Without loss of generality we can assume $\tilde{\Delta}_\nu \geq 4$. Put $\upsilon = 2\lfloor \tilde{\Delta}_\nu/2 \rfloor \geq \tilde{\Delta}_\nu/2$, Then $\{x: 0 \leq x < c_0\sqrt{\upsilon}\} \supseteq \{x: 0 \leq x < c_0\sqrt{\tilde{\Delta}_\nu/2}\}$. Therefore, from Lemma 2.1 of SS(1991) it follows

$$P\{\xi > x\} = (1-\Phi(x))\exp\{L(x)\}\left(1+\theta_1 c_2(x+1)\Delta^{-1/(1+2\nu)}\right), \text{ for } 0 \leq x < c_1\Delta^{1/(1+2\nu)},$$

where $c_1 = c_0/\sqrt{2 \cdot 18^{1/(1+2\nu)}}$ and $c_2 = \sqrt{2 \cdot 18^{1/(1+2\nu)}}\,f_1(x,(\Delta^2/18)^{1/(1+2\nu)})$. At last, exactly as in SS (1991, p. 35-37) we have $L(x) = L_\nu(x) + \theta(x/\Delta_\nu)^3$. Assertion 2 follows.

The following relation is well-known, see SS (1991, Eq. (1.34)).

**Assertion 3.** One has

$$\mathcal{C}_k(\varsigma) = k!\sum (-1)^{m_1+m_2+...+m_k-1}(m_1+m_2+...+m_k-1)!\prod_{l=1}^{k}\frac{1}{m_l!}\left(\frac{\alpha_l(\varsigma)}{l!}\right)^{m_l}$$

here $\Sigma$ is summation over all non-negative integer $m_1, m_2, ..., m_k$ such that $m_1 + 2m_2 + ... + km_k = k$.